# OPTIMAL CONTROL OF VOLTERRA EQUATIONS WITH IMPULSES

S. A. Belbas  
Mathematics Department  
University of Alabama  
Box 870350  
Tuscaloosa, AL. 35487-0350.  
USA.

W. H. Schmidt  
Institut für Mathematik und Informatik  
Universität Greifswald  
Jahnstraße 15a  
D-17487 Greifswald.  
GERMANY.

e-mail:  
SBELBAS@GP.AS.UA.EDU,  
SBELBAS@BAMA.UA.EDU

Abstract

We consider an optimal control problem for a system governed by a Volterra integral equation with impulsive terms. The impulses act on both the state and the control; the control consists of switchings at discrete times. The cost functional includes both, an integrated cost rate (continuous part) and switching costs at the discrete impulse times (discrete part). We prove necessary optimality conditions of a form analogous to a discrete maximum principle. For the particular case of a system governed by impulsive ordinary differential equations, we obtain an impulsive maximum principle as a special case of the necessary optimality conditions for impulsive Volterra equations.

1. Introduction.

In this paper, we obtain necessary conditions for optimality for Volterra integral equations with impulsive terms and piecewise constant controls. The cost functional that is to be minimized includes an integral term for the rate of cost per unit time, as well as discrete switching costs at the times of changing the values of the control. For simplicity, we consider only discrete controls and fixed switching times. The combination of continuous and discrete controls, and the case of variable switching times, are objects of ongoing work and will be reported elsewhere.

The novelty of the problem we consider here lies in the inclusion of impulsive terms in the Volterra integral equation. Without impulsive terms in the state dynamics, the case of piecewise constant controls has been treated in [S1, S2, S3]. We note here that the general questions of optimality conditions for optimal control problems governed by Volterra integral equations have been studied in, among other references, [B, BA, CA, RS, S1, S2, S3, S4].

It is well known that Volterra integral equations can be used to model many classes of phenomena, for example population dynamics, continuum mechanics of materials with memory, economic problems, the spread of epidemics, non-local problems of diffusion and heat conduction, etc. Some of these applications may be found in the classical references [D, K], and other examples can be found in more specialized works. The corresponding control problems for such systems lead to optimal control problems for Volterra integral equations. An explicit example of applying the methods of optimal control of Volterra integral equations to economics is the paper [KM].

All kinds of control systems can be subjected to impulsive conditions. Impulsive control systems may arise from state models that are intrinsically impulsive, i.e. the physical model without a control function (the uncontrolled model) still involves impulsive terms, as it happens, for example, in systems of variable structure when a system involves transitions through different operating regimes. Impulsive control problems can also arise because the action of the control may involve impulses applied to the state of the system, as, for example, investment decisions in economics, or the injection of a medical drug into a patient in mathematical models in pharmacology. Therefore, it is of interest to study the optimal control of impulsive Volterra integral equations.

The mathematical analysis of problems with discrete controls is useful not only when the original formulation of the state dynamics of a controlled system involves discrete controls, but also when one tries to numerically solve an optimal control problem that was originally formulated as a problem with continuous controls. The numerical solution of optimal control problems often involves a discretization of the control function (see, e.g., [VS] for discrete approximations to ODE control problems), and in that case one has to solve the corresponding optimal control problem with a discrete-valued control function.

The problem of impulsive control of a Volterra integral equation with fixed instants of impulses or switching of the control leads to a novel type of necessary conditions of optimality, analogous to the discrete version of Pontryagin's maximum principle. The discrete maximum principle, for systems governed by finite-difference equations (the discrete analogue of ordinary differential equations) has been treated by several researchers; rigorous results have been obtained notably in [BI, BO]. The problem we consider in this paper also lead to a discrete version of a maximum principle, but the state dynamics in our problem is a Volterra equation in continuous time with additional impulsive terms. To the best of our knowledge, this problem has not been previously considered in the researcg literature. The results of this paper rely on a thorough analysis of the associated variational equations for impulsive Volterra integral equations. These results are proved in sections 4 through 6 of this paper. For the particular case of impulsive ordinary differential equations, we prove, in section 7, a version of a maximum principle for impulsive ODE as a particular case of the results for impulsive Volterra integral equations.

## 2. Statement of the problem.

We consider a controlled Volterra equation with impulses at specified times $0 < \tau_1 < \tau_2 < \cdots < \tau_N < T$, where [0, T] is the time-horizon of the optimal control problem. The system dynamics is described by

$$x(s) = h(s) + \int_0^s f(s,t,x(t),u(t))dt + \sum_{i: 0 < \tau_i < s} G(s,\tau_i, x(\tau_i^-), u(\tau_i^-), u(\tau_i^+))$$

(2.1)

where the control function u(.) is piecewise constant,

$$u(t) = a_i \text{ for } \tau_{i-1} \le t < \tau_i$$

(2.2)

Of course, $u(\tau_i^-) = a_i$, $u(\tau_i^+) = a_{i+1}$, for $i = 1, 2, ..., N$. We set $\tau_0 := 0$, $\tau_{N+1} := T$.

The control function u takes values in a compact set $U \subset \mathrm{IR}^m$. A control policy is identified with an ordered (N+1)-tuple of elements of U, $(a_1, a_2, \cdots, a_N, a_{N+1})$ so that the associated control function $u(\cdot)$ is given by (2.2) for i=1,2,...,N+1.

The functions f and G satisfy the following conditions:

(C1) f(s,t,x,a) and G(s,t,x,a,b) are continuous for $0 \le t \le s \le T$, x in $\mathrm{IR}^n$, a and b in U.

(C2) f and G are Lipschitz in x, uniformly in s, t, a, b :

$$|f(s,t,x_1,\alpha) - f(s,t,x_2,\alpha)| \le C_f |x_1 - x_2|$$
$$|G(s,t,x_1,\alpha,\beta) - G(s,t,x_2,\alpha,\beta)| \le C_G |x_1 - x_2|$$

(2.3)

The objective of the optimal control problem is to minimize a functional J defined by

$$J := \int_0^T g(t, x(t), u(t))dt + \varphi(x(\tau_1^-), x(\tau_2^-), ..., x(\tau_{N+1}^-), a_1, a_2, ..., a_{N+1})$$

(2.4)

## 3. Existence and elementary properties of solutions of the state equation.

We use the definitions and notation of the previous section. We set
$\boldsymbol{\tau} = (\tau_1, \tau_2, ..., \tau_{N+1})$,
$\mathbf{a} = (a_1, ..., a_{N+1})$.
For each collection of impulse times and corresponding value of the control, we seek a solution of (2.1) in the space $C(0,T;\boldsymbol{\tau})$ of functions $x(t)$ that are bounded on $[0,T]$, continuous on each interval $(\tau_{i-1}, \tau_i)$, $i = 1, 2, ..., N+1$, and have limits
$x(\tau_0^+)$, $x(\tau_{N+1}^-)$, $x(\tau_i^\pm)$ for $1 \leq i \leq N$ (where $x(\tau^\pm) := \lim_{t \to \tau^\pm} x(t)$).

**Theorem 3.1.** Under conditions (C1, C2) of section 2, for every $\boldsymbol{\tau}$ and $\mathbf{a}$, equation (2.1) possesses a unique solution $x(.)$ in the space $C(0,T;\boldsymbol{\tau})$.

Proof: Define $y_i \in C([\tau_{i-1}, \tau_i] \mapsto \mathbb{R})$ as follows:

$y_1$ is the unique solution (guaranteed by conditions C1, C2) of the Volterra equation

$$y_1(s) = h(s) + \int_0^s f(t, s, x(t), a_1) dt, \quad s \in [0, \tau_1] \tag{3.1}$$

and, inductively, for $i = 1, 2, ..., N$, $y_{i+1}$ is solution of

$$y_{i+1}(s) = h(s) + \sum_{j=1}^{i} \left[ \int_{\tau_{j-1}}^{\tau_j} f(s, t, y_j(t), a_j) dt + G(s, \tau_j, y_j(\tau_j), a_j, a_{j+1}) \right] +$$
$$+ \int_{\tau_i}^{s} f(s, t, y_{i+1}(t), a_{j+1}) dt, \quad \text{for } s \in [\tau_i, \tau_{i+1}] \tag{3.2}$$

Then the function $x(s)$ defined by

$$x(s) = y_i(s) \text{ for } s \in (\tau_{i-1}, \tau_i); \; i = 1, 2, ..., N+1;$$
$$x(\tau_i^-) = y_i(\tau_i); \; x(\tau_i^+) = y_i(\tau_i) + G(\tau_i, \tau_i, y_i(\tau_i), a_i, a_{i+1}) \tag{3.3}$$

solves the impulsive Volterra equation (2.1). ///

For the purpose of devising iterative methods for solving impulsive Volterra equations, it is useful to obtain the solution of an impulsive Volterra equation via an iterative method. First, we notice that, under the stated assumptions, we have

$$x(\tau_i^-) = h(\tau_i) + \int_0^{\tau_i} K(\tau_i, t, x(t), u(t))dt + \sum_{j=1}^{i-1} G(\tau_i, \tau_j, x(\tau_j^-), u(\tau_j^-), u(\tau_j^+))$$

We consider the space $V := C(0,T;\boldsymbol{\tau}) \times \mathrm{IR}^N$ with a vector-valued norm, defined for each $(x, \eta)$ in V (with $x \in C(0,T;\boldsymbol{\tau})$, $\eta = (\eta_1, \eta_2, ..., \eta_N) \in \mathrm{IR}^N$) by

$$\|(\xi, \eta)\|_\mu = \begin{bmatrix} \|\xi\|_\mu \\ \|\eta\|_\mu \end{bmatrix}; \quad \|\xi\|_\mu := \sup_{0 \le t \le T} \left( e^{-\mu t} |\xi(t)| \right); \quad \|\eta\|_\mu = \max_{1 \le i \le N} \left( e^{-\mu \tau_i} |\eta_i| \right)$$

Even though, for simplicity, we use the same symbol for all 3 norms, the norm $\|\cdot\|_\mu$ on $C(0,T;\boldsymbol{\tau})$, the norm $\|\cdot\|_\mu$ on $\mathrm{IR}^N$, and the $\mathrm{IR}^2$-valued norm $\|\cdot\|_\mu$ on $C(0,T;\boldsymbol{\tau}) \times \mathrm{IR}^N$, it should be clear from the context which norm is used in each case.

The problem of solving the impulsive Volterra equation (2.1) amounts to finding a fixed point of an operator S defined on V by

$$S(\xi, \eta) := \begin{bmatrix} S_c(\xi, \eta) \\ S_d(\xi, \eta) \end{bmatrix};$$

$S_c : V \mapsto C(0, T; \boldsymbol{\tau})$ ; $S_d : V \mapsto \mathrm{IR}^N$ ;

$$(S_c(\xi, \eta))(s) := h(s) + \int_0^s f(s, t, \xi(t), u(t))dt + \sum_{j:0 < \tau_j < s} G(s, \tau_j, \eta_j, u(\tau_j^-), u(\tau_j^+));$$

$$(S_d(\xi, \eta))_i = h(\tau_i) + \int_0^{\tau_i} f(\tau_i, t, \xi(t), u(t))dt + \sum_{j:1 \le j \le i-1} G(\tau_i, \tau_j, \eta_j, u(\tau_j^-), u(\tau_j^+))$$

We follow the standard convention that a sum over an empty set of indices is, by definition, equal to zero; for example, the sum that appears in the definition of $S_d$ above is zero when i=1.

We shall call $S_c$ and $S_d$ the continuous and the discrete components of S.

We have:

Lemma 3.1. Solving (3.1) is equivalent to finding a fixed point of the operator S.

<u>Proof:</u> If x is a solution of (3.1), then we can take $\xi^*(s) = x(s)$ for $s \notin \tau$ and $\eta_i^* = x(\tau_i^-)$ for $i = 1, 2, ..., N$; then $(\xi^*, \eta^*)$ is a fixed point of S.

Conversely, if $(\xi^*, \eta^*)$ is a fixed point of S in V, then, for $s \in [0, \tau_1)$, we have

$$\xi^*(s) = h(s) + \int_0^s f(s, t, \xi^*(t), u(t)) dt$$

and at $\tau_1$ we have

$$\eta_1^* = h(\tau_1) + \int_0^{\tau_1} f(\tau_1, t, \xi^*(t), u(t)) dt$$

therefore $\eta_1^* = \xi^*(\tau_1^-)$. Inductively, if

$$\xi^*(s) = h(s) + \int_0^s f(s, t, \xi^*(t), u(t)) dt + \sum_{j: 0 < \tau_j < s} G(s, \tau_j, \xi^*(\tau_j^-), u(\tau_j^-), u(\tau_j^+))$$

for all $s < \tau_k$, $s \notin \tau$, and $\eta_j^* = \xi^*(\tau_j^-)$ for $j = 1, 2, ..., k-1$, then we have, at $\tau_k$,

$$\eta_k^* = h(\tau_k) + \int_0^{\tau_k} f(\tau_k, t, \xi^*(t), u(t)) dt +$$
$$+ \sum_{j=1}^{k-1} G(\tau_k, \tau_j, \eta_j^*, u(\tau_j^-), u(\tau_j^+)) =$$
$$= h(\tau_k) + \int_0^{\tau_k} f(\tau_k, t, \xi^*(t), u(t)) dt +$$
$$+ \sum_{j=1}^{k-1} G(\tau_k, \tau_j, \xi^*(\tau_j^-), u(\tau_j^-), u(\tau_j^+)) = \xi^*(\tau_k^-)$$

and, for $s \in (\tau_k, \tau_{k+1})$,

$$\xi^*(s) = h(s) + \int_0^s f(s, t, \xi^*(t), u(t)) dt + \sum_{j=1}^k G(s, \tau_j, \eta_j^*, u(\tau_j^-), u(\tau_j^+)) =$$
$$= h(s) + \int_0^s f(s, t, \xi^*(t), u(t)) dt + \sum_{j=1}^k G(s, \tau_j, \xi^*(\tau_j^-), u(\tau_j^-), u(\tau_j^+))$$

so that $\xi^*$ solves the impulsive Volterra equation on $[0, \tau_{k+1})$. The induction is complete, thus showing that $\xi^*$ solves the impulsive Volterra equation on $[0, T]$ and $\eta_j^* = \xi^*(\tau_j^-)$ for all j=1, 2, ..., N. ///

Lemma 3.2. If $\mu>0$ is sufficiently large, then, for any two elements $(\xi^1, \eta^1)$, $(\xi^2, \eta^2)$ of V, we have

$$\|S_c(\xi^1, \eta^1) - S_c(\xi^2, \eta^2)\|_\mu \leq a_{11} \|\xi^1 - \xi^2\|_\mu + a_{12} \|\eta^1 - \eta^2\|_\mu ;$$
$$\|S_d(\xi^1, \eta^1) - S_d(\xi^2, \eta^2)\|_\mu \leq a_{21} \|\xi^1 - \xi^2\|_\mu + a_{22} \|\eta^1 - \eta^2\|_\mu$$

(3.4)

where all constants $a_{ij}$ are nonnegative and the eigenvalues $\lambda_1$, $\lambda_2$ of the matrix

$$A := \begin{bmatrix} a_{11} & a_{12} \\ a_{21} & a_{22} \end{bmatrix} \text{ satisfy } |\lambda_i| < 1 \text{ for i=1, 2.}$$

Note: The nonnegativity of all elements of A is equivalent to A being an order-preserving operator on $\mathrm{I\!R}^2$ with respect to the coordinate-wise partial order in $\mathrm{I\!R}^2$.

Proof: We recall (cf., e.g., [DM]) that the necessary and sufficient condition for a real $2\times 2$ matrix A to have eigenvalues with moduli less than 1 is $|\mathrm{tr}(A)| - 1 < \det(A) < 1$, where $\mathrm{tr}(A)$ and $\det(A)$ are the trace and the determinant of A.

Let $h := \min_{1 \leq i \leq N+1} (\tau_i - \tau_{i-1})$, and, for each s, let $N_s := \max\{i \geq 0 : \tau_i < s\}$. Then we have

$$e^{-\mu s} |(S_c(\xi_1, \eta_1))(s) - (S_c(\xi_2, \eta_2))(s)| \leq$$

$$\leq C_f e^{-\mu s} \int_0^s |\xi_1(t) - \xi_2(t)| dt + C_G e^{-\mu s} \sum_{i=1}^{N_s} |\eta_{1,i} - \eta_{2,i}| \leq$$

$$\leq C_f e^{-\mu s} \int_0^s e^{\mu t} \|\xi_1 - \xi_2\|_\mu dt + C_G e^{-\mu s} \sum_{i=1}^{N_s} e^{\mu \tau_i} \|\eta_{1,i} - \eta_{2,i}\|_\mu$$

Now,

$$e^{-\mu s} \int_0^s e^{\mu t} \|\xi_1 - \xi_2\|_\mu dt = \frac{1 - e^{-\mu s}}{\mu} \|\xi_1 - \xi_2\|_\mu \leq \frac{1 - e^{-\mu T}}{\mu} \|\xi_1 - \xi_2\|_\mu$$

and

$$e^{-\mu s}\sum_{i=1}^{N_s} e^{\mu\tau_i} \|\eta_1 - \eta_2\|_\mu = \left[\sum_{i=1}^{N_s} e^{-\mu(s-\tau_i)}\right]\|\eta_1 - \eta_2\|_\mu \leq$$

$$\leq \left[1 + \sum_{j=1}^{N_s-1} e^{-\mu j h}\right]\|\eta_1 - \eta_2\|_\mu = \left[1 + e^{-\mu h}\frac{1-e^{-(N_s-1)\mu h}}{1-e^{-\mu h}}\right]\|\eta_1 - \eta_2\|_\mu \leq$$

$$\leq \left[1 + e^{-\mu h}\frac{1-e^{-(N-1)\mu h}}{1-e^{-\mu h}}\right]\|\eta_1 - \eta_2\|_\mu$$

thus

$$\|S_c(\xi_1,\eta_1) - S_c(\xi_2,\eta_2)\|_\mu \leq C_f \frac{1-e^{-\mu T}}{\mu}\|\xi_1-\xi_2\|_\mu +$$

$$+ \left[1 + e^{-\mu h}\frac{1-e^{-(N-1)\mu h}}{1-e^{-\mu h}}\right]\|\eta_1-\eta_2\|_\mu$$

Also, by the same type of calculation,

$$e^{-\mu\tau_i}|(S_c(\xi_1,\eta_1))_i - (S_c(\xi_2,\eta_2))_i| \leq C_f \frac{1-e^{-\mu\tau_i}}{\mu}\|\xi_1-\xi_2\|_\mu +$$

$$+ C_G \sum_{j=1}^{i-1} e^{-\mu(\tau_i-\tau_j)}\|\eta_1-\eta_2\|_\mu$$

and

$$\sum_{j=1}^{i-1} e^{-\mu(\tau_i-\tau_j)} \leq \sum_{j=1}^{i-1} e^{-\mu j h} \leq \sum_{j=1}^{N-1} e^{-\mu j h} = e^{-\mu h}\frac{1-e^{-\mu(N-1)h}}{1-e^{-\mu h}}$$

thus

$$\|S_c(\xi_1,\eta_1) - S_c(\xi_2,\eta_2)\|_\mu \leq C_f \frac{1-e^{-\mu\tau_N}}{\mu}\|\xi_1-\xi_2\|_\mu +$$

$$+ C_G e^{-\mu h}\frac{1-e^{-\mu(N-1)h}}{1-e^{-\mu h}}\|\eta_1-\eta_2\|_\mu$$

Thus, we can take, in the inequalities (3.4),

$$a_{11} \equiv a_{11}(\mu) = C_f \frac{1-e^{-\mu T}}{\mu}; \ a_{12} \equiv a_{12}(\mu) = \left[1 + e^{-\mu h} \frac{1-e^{-(N-1)\mu h}}{1-e^{-\mu h}}\right];$$

$$a_{21} \equiv a_{21}(\mu) = C_f \frac{1-e^{-\mu \tau_N}}{\mu}; \ a_{22} \equiv a_{22}(\mu) = C_G e^{-\mu h} \frac{1-e^{-\mu(N-1)h}}{1-e^{-\mu h}}$$

(3.5)

It follows from (3.5) above that $(a_{11}(\mu), a_{21}(\mu), a_{22}(\mu)) \to (0,0,0)$ as $\mu \to \infty$, and $a_{12}(\mu)$ remains bounded as $\mu \to \infty$; consequently, $(\text{tr}(A), \det(A)) \to (0, 0)$ as $\mu \to \infty$, and therefore, for $\mu$ sufficiently large, the inequalities $|\text{tr}(A)| - 1 < \det(A) < 1$ will be satisfied. ///

We consider the following iterative scheme for the solution of (2.1):

$(\xi_{(0)}, \eta_{(0)})$ is an arbitrary element of V;

for k=0, 1, 2, ..., $(\xi_{(k+1)}, \eta_{(k+1)})$ is defined by

$$\xi_{(k+1)}(s) = h(s) + \int_0^s f(s,t,\xi_{(k)}(t),u(t))dt + \sum_{i:0<\tau_i<s} G(s,\tau_i,\eta_{(k),i},u(\tau_i^-),u(\tau_i^+));$$

$$\eta_{(k+1),i} = h(\tau_i) + \int_0^{\tau_i} f(\tau_i,t,\xi_{(k)}(t),u(t))dt + \sum_{j=1}^{i-1} G(\tau_i,\tau_j,\eta_{(k),j},u(\tau_j^-),u(\tau_j^+))$$

(3.6)

Then we have:

Theorem 3.2. Under the conditions (C1, C2), the iterative method defined by (3.6) above converges to $(x(\cdot),(x(\tau_1^-), x(\tau_2^-),..., x(\tau_N^-)))$, as $k \to \infty$; the convergence of $\xi_{(k)}$ to $x(\cdot)$ is uniform on each interval $[\tau_{i-1}, \tau_i]$, $i = 1,2,..., N+1$, in the sense that, if we define the restrictions $\xi_{(k)}^i, x^i$ of $\xi_{(k)}, x(\cdot)$ to the interval $[\tau_{i-1}, \tau_i]$ by

$\xi_{(k)}^i(t) = \xi_{(k)}(t), x^i(t) = x(t),$ for $t \in (\tau_{i-1}, \tau_i)$;

$\xi_{(k)}^i(\tau_{i-1}) = \xi_{(k)}(\tau_{i-1}^+); x^i(\tau_{i-1}) = x(\tau_{i-1}^+); \xi_{(k)}^i(\tau_i) = \xi_{(k)}(\tau_i^-); x^i(\tau_i) = x(\tau_i^-)$

then $\xi_{(k)}^i \to x^i$ as $k \to \infty$ uniformly on $[\tau_{i-1}, \tau_i]$ for all i=1, 2, ..., N+1.

Proof: By lemma 3.2, if $\mu$ is sufficiently large, the operator S is a contraction with respect to the vector-valued norm $\|\cdot\|_\mu$ on V. The contraction property with respect to the vector-valued norm $\|\cdot\|_\mu$ means that, for all z, w in V, we have $\|Sz - Sw\|_\mu \leq A \|z - w\|_\mu$ where the 2-by-2 real matrix A has nonnegative elements and the eigenvalues of A have moduli less than 1. Consequently, the iterates of S, with arbitrary initial data, converge to the unique fixed point of S in the topology induced on V by the vector-valued norm $\|\cdot\|_\mu$; this is a well-known extension of the standard Banach fixed-point theorem to the case of a vector-valued metric, and the proof proceeds as in the standard case. Convergence with respect to the norm $\|\cdot\|_\mu$ on C(0,T;$\boldsymbol{\tau}$) is equivalent to uniform convergence on each $[\tau_{i-1}, \tau_i]$. The fixed point of S gives the solution of (3.1) by lemma 3.1.  ///

Remark 3.1.  It follows from (3.6) that $\eta_{(k+1),i} = \xi_{(k+1)}(\tau_i^-)$ for all k$\geq$0, so that, for k$\geq$1, (3.6) can also be written in the form

$$\xi_{(k+1)}(s) = h(s) + \int_0^s f(s,t,\xi_{(k)}(t),u(t))dt + \sum_{i:0<\tau_i<s} G(s,\tau_i,\xi_{(k)}(\tau_i^-),u(\tau_i^-),u(\tau_i^+))$$

Of course, if $\eta_{(0)}$ is chosen as $\eta_{(0),i} = \xi_{(0)}(\tau_i^-)$, then (3.) holds for all k=0, 1, 2, ... . ///

4. Linear Volterra equations with impulses.

The theory of ordinary Volterra equations (without impulses) contains particular results for the case of linear integral equations, including the method of iterated kernels and the duality properties of the resolvent. Inasmuch as the variational equations associated with an impulsive optimal control problem are linear equations, we have to establish the corresponding results for linear impulsive Volterra equations.

A linear Volterra equation with impulses has the form

$$x(s) = h(s) + \int_0^s K(s,t)x(t)dt + \sum_{i: 0 < \tau_i < s} L(s,\tau_i)x(\tau_i^-)$$

(4.1)

where

$$x(\tau_i^-) = h(\tau_i) + \int_0^{\tau_i} K(\tau_i,t)x(t)dt + \sum_{j: j < i} L(\tau_i,\tau_j)x(\tau_j^-)$$

(4.2)

We can represent the solution of (4.1) by using at first convolutions of the impulsive parts of the operator M, i.e. the terms containing the functions $L(s,\tau_i)$ ; this leads to a discrete resolvent. For every pair of indices (i, j) with i>j, we define the set **P**(j, i) of increasing paths from j to i as the set of all ordered collections s of the form

$$\sigma := (\tau_j, \tau_{k_1}, \tau_{k_2}, ..., \tau_{k_\alpha}, \tau_i) ; j < k_1 < k_2 < ... < k_\alpha < i$$

and we define the discrete resolvent kernel L by

$$\Lambda(\sigma) := L(\tau_i, \tau_{k_\alpha})L(\tau_{k_\alpha}, \tau_{k_{\alpha-1}})...L(\tau_{k_1}, \tau_j)$$

We extend the definition of L($\sigma$) to the case j=i by setting L(s)=1 if $\sigma$={ $\tau_i$ }.
Then we can show:

Theorem 4.1. Equation (4.2) is equivalent to

$$x(\tau_i^-) = \sum_{j=1}^{i} \sum_{\sigma \in \mathbf{P}(j,i)} \Lambda(\sigma)[h(\tau_j) + \int_0^{\tau_j} K(\tau_j,t)x(t)dt] \ , \ i=1, 2, ..., N$$

(4.3)

i.e., if x(.) is the solution of (4.1) then the quantities $x(\tau_i^-)$ can be calculated either from (4.2) or from (4.3).

Proof: For i=1, we have, from (4.1),

$$x(\tau_1^-) = h(\tau_1^-) + \int_0^{\tau_1} K(\tau_1,t)x(t)dt$$

which agrees with (4.3) for i=1. Inductively, if

$$x(\tau_j^-) = \sum_{k=1}^{j} \sum_{\sigma \in \mathbf{P}(k,j)} \Lambda(\sigma)[h(\tau_k) + \int_0^{\tau_k} K(\tau_k,t)x(t)dt]$$
$$\forall j = 1,2,...,i-1$$

then we substitute the above into (4.2) and we get

$$x(\tau_i^-) = h(\tau_i) + \int_0^{\tau_i} K(\tau_i,t)x(t)dt + \sum_{j:j<i} \sum_{k=1}^{j} \sum_{\sigma \in \mathbf{P}(k,j)} L(\tau_i,\tau_j)\Lambda(\sigma) \cdot$$
$$\cdot [h(\tau_k) + \int_0^{\tau_k} K(\tau_k,t)x(t)dt]$$

Now, every σ' in $\mathbf{P}$(k,i) with k<i is obtained as concatenation of some σ in $\mathbf{P}$(k,j) for some j<i (including the possibility j=k) and the additional $\tau_i$ after the last $\tau_j$ of σ. Therefore, for such a σ', we have

$$\Lambda(\sigma') = L(\tau_i,\tau_j)\Lambda(\sigma)$$

For k=i, the only element of $\mathbf{P}$(k, i) is {$\tau_i$}.
By making a change in the order of summation, we have

$$\sum_{j:j<i} \sum_{k=1}^{j} \sum_{\sigma \in \mathbf{P}(k,j)} L(\tau_i, \tau_j) \Lambda(\sigma) = \sum_{k=1}^{i} \sum_{\sigma' \in \mathbf{P}(k,i)} \Lambda(\sigma')$$

and consequently

$$x(\tau_i^-) = h(\tau_i) + \int_0^{\tau_i} K(\tau_i, t) x(t) dt +$$

$$+ \sum_{k=1}^{i} \sum_{\sigma' \in \mathbf{P}(k,i)} \Lambda(\sigma')[h(\tau_k) + \int_0^{\tau_k} K(\tau_k, t) x(t) dt]$$

The induction is complete. ///

Remark 4.1.  An expression like

$$\sum_{j=1}^{i} \sum_{\sigma \in \mathbf{P}(j,i)} \Lambda(\sigma) \xi_j$$

for any variables $\xi_j : j = 1, 2, ..., N$, can also be written as

$$\xi_i + \sum_{m=1}^{\infty} \sum_{(k_1, k_2, ..., k_m)} \sum_{(\tau_{k_1}, \tau_{k_2}, \cdots, \tau_{k_m})} [\xi_{k_1} L(\tau_{k_1}, \tau_{k_2}) L(\tau_{k_2}, \tau_{k_3}) \cdots$$
$$\cdots L(\tau_{k_{m-1}}, \tau_{k_m}) L(\tau_{k_m}, \tau_i)]$$

(4.4)

since $L(\tau_k, \tau_i) = 0$ if $k \geq i$. The third summation in (4.4) above is taken over all $(\tau_{k_1}, \tau_{k_2}, \cdots, \tau_{k_m})$ in $\{\tau_1, \tau_2, ..., \tau_N\}^m$ . ///

As a consequence of theorem 4.2, equation (4.1) takes the form

$$x(s) = h(s) + \int_0^s K(s,t)x(t)dt + \sum_{i:0<\tau_i<s} \sum_{j=1}^{i} \sum_{\sigma \in \mathbf{P}(j,i)} L(s,\tau_i)\Lambda(\sigma)[h(\tau_j) +$$
$$+ \int_0^{\tau_j} K(\tau_j,t)x(t)dt]$$
(4.5)

The calculations can be somewhat simplified if we make the following conventions: we extend K(s,t) so that it is defined for all $(s,t) \in [0,T] \times [0,T]$ with $K(s,t) = 0$ for $t \geq s$, and similarly we extend $L(s,\tau_i)$ for all i=1,2,...,N with the convention that $L(s,\tau_i) = 0$ if $\tau_i \geq s$.

By using these extensions of K and L, we can replace the integral $\int_0^{\tau_j} K(\tau_j,t)x(t)dt$ on the right-hand side of (4.2) above by $\int_0^s K(\tau_j,t)x(t)dt$ without violating the validity of (4.2). Then (4.5) is a linear Volterra integral equation of the form

$$x(s) = \mathbf{h}(s) + \int_0^s \mathbf{K}(s,t)x(t)dt$$
(4.6)

with discontinuous kernel K given by

$$\mathbf{K}(s,t) := K(s,t) + \sum_{i:0<\tau_i<s} \sum_{j=1}^{i} \sum_{\sigma \in \mathbf{P}(j,i)} L(s,\tau_i)\Lambda(\sigma)K(\tau_j,t)$$
(4.7)

and discontinuous forcing term

$$\mathbf{h}(s) := h(s) + \sum_{i:0<\tau_i<s} \sum_{j=1}^{i} \sum_{\sigma \in \mathbf{P}(j,i)} L(s,\tau_i)\Lambda(\sigma)h(\tau_j)$$
(4.8)

We extend the definition of $\mathbf{P}(j, i)$ to all pairs (j, i) by setting $\mathbf{P}(j, i) := \emptyset$ when i<j.

We note that, in expressions like $\sum_{i:0<\tau_i<s} \sum_{j=1}^{i} \sum_{\sigma \in \mathbf{P}(j,i)} L(s,\tau_i)\Lambda(\sigma)K(\tau_j,t)$ or $\sum_{i:0<\tau_i<s} \sum_{j=1}^{i} \sum_{\sigma\in\mathbf{P}(j,i)} L(s,\tau_i)\Lambda(\sigma)h(\tau_j)$ , the expression inside the summation becomes zero when i<j or s≥ $\tau_i$ . Then we have

$$\sum_{i:0<\tau_i<s} \sum_{j=1}^{i} \sum_{\sigma\in\mathbf{P}(j,i)} L(s,\tau_i)\Lambda(\sigma)K(\tau_j,t) = \sum_{i=1}^{N} \sum_{j=1}^{N} \sum_{\sigma\in\mathbf{P}(j,i)} L(s,\tau_i)\Lambda(\sigma)K(\tau_j,t);$$

$$\sum_{i:0<\tau_i<s} \sum_{j=1}^{i} \sum_{\sigma\in\mathbf{P}(j,i)} L(s,\tau_i)\Lambda(\sigma)h(\tau_j) = \sum_{i=1}^{N} \sum_{j=1}^{N} \sum_{\sigma\in\mathbf{P}(j,i)} L(s,\tau_i)\Lambda(\sigma)h(\tau_j)$$

Further, we make the convention to write $\sum_{i,j,\sigma}$ for $\sum_{i=1}^{N} \sum_{j=1}^{N} \sum_{\sigma\in\mathbf{P}(j,i)}$ .

We shall prove:

Theorem 4.2. The problem (4.1, 4.2) is equivalent to the problem (4.6, 4.7, 4.8).

Proof: By theorem 4.1, if x(.) solves (4.2, 4.3), then it also solves (4.6, 4.7, 4.8); it remains to show that, if x(.) solves (4.6, 4.7, 4.8), and if each $x(\tau_i^-)$ is defined by

$$x(\tau_i^-) = \mathbf{h}(\tau_i^-) + \int_0^s \mathbf{K}(\tau_i^-,t)x(t)dt$$

(4.9)

then x(.) must also solve

$$x(s) = h(s) + \int_0^s K(s,t)x(t)dt + \sum_{i:0<\tau_i<s} L(s,\tau_i)x(\tau_i^-)$$

(4.10)

We assume that x(.) solves (4.6, 4.7, 4.8). We have

$$x(\tau_k^-) = \mathbf{h}(\tau_k^-) + \int_0^{\tau_k} \mathbf{K}(\tau_k^-, t)x(t)dt =$$

$$= h(\tau_k) + \sum_{i,j,\sigma} L(\tau_k, \tau_i)\Lambda(\sigma)h(\tau_j) +$$

$$+ \int_0^{\tau_k} [K(\tau_k, t) + \sum_{i,j,\sigma} L(\tau_k, \tau_i)\Lambda(\sigma)K(\tau_j, t)]x(t)dt$$

(4.11)

and consequently

$$L(s,\tau_k)x(\tau_k^-) = L(s,\tau_k)h(\tau_k) + \sum_{i,j,\sigma} L(s,\tau_k)L(\tau_k,\tau_i)\Lambda(\sigma)h(\tau_j) +$$

$$+ \int_0^s [L(s,\tau_k)K(\tau_k,t) + \sum_{i,j,\sigma} L(s,\tau_k)L(\tau_k,\tau_i)\Lambda(\sigma)K(\tau_j,t)]x(t)dt$$

(4.12)

Since $x(.)$ solves (4.6), we have

$$x(s) = h(s) + \sum_{i,j,\sigma} L(s,\tau_i)\Lambda(\sigma)h(\tau_j) +$$

$$+ \int_0^s [K(s,t) + \sum_{i,j,\sigma} L(s,\tau_i)\Lambda(\sigma)K(\tau_j,t)]x(t)dt$$

(4.13)

From (4.12), and also taking into account the properties of $\Lambda(\sigma)$ that were established in the proof of theorem 4.1, we have

$$\sum_k L(s,\tau_k)x(\tau_k^-) =$$

$$= \sum_k [L(s,\tau_k)h(\tau_k) + \sum_{i,j,\sigma} L(s,\tau_k)L(\tau_k,\tau_i)\Lambda(\sigma)h(\tau_j) +$$

$$+ \int_0^s \{L(s,\tau_k)K(\tau_k,t) + \sum_{i,j,\sigma} L(s,\tau_k)L(\tau_k,\tau_i)\Lambda(\sigma)K(\tau_j,t)\}x(t)dt] =$$

$$= \sum_{i',j,\sigma'} L(s,\tau_{i'})\Lambda(\sigma')h(\tau_j) +$$

$$+ \int_0^s [K(s,t) + \sum_{i',j,\sigma'} L(s,\tau_{i'})\Lambda(\sigma')K(\tau_j,t)]x(t)dt$$

(4.14)

It follows from (4.13) and (4.14) that x(.) solves (4.1). The fact that $x(\tau_k^-)$, as evaluated in (4.11), coincides with the $x(\tau_k^-)$ that would be evaluated according to (4.2) follows by a similar argument.  ///

The solution of (4.6, 4.7, 4.8) can be expressed by means of a resolvent **R**(s, t) defined via iterated kernels computed from the discontinuous kernel **K**(s, t). We define the convolution of two kernels $\mathbf{K}_1(.,.)$ and $\mathbf{K}_2(.,.)$ by

$$(\mathbf{K}_1 * \mathbf{K}_2)(s,t) = \int_t^s \mathbf{K}_1(s,t_1)\mathbf{K}_2(t_1,t)dt_1$$

(4.15)

The class of kernels for which the convolution (4.15) is defined consists of functions **K**(s, t) that are bounded on 0<t<s<T, are continuous in the second variable, and have at most jump discontinuities at the points $s = \tau_i$, $i = 1,2,...,N$. It is readily verified that this convolution is associative, and therefore it is meaningful to define the convolutional powers $\mathbf{K}^{*n}$, where the exponent signifies n-fold convolution of **K** with itself. We have:

<u>Theorem 4.3.</u>  Define the resolvent kernel $\mathbf{R}(s,t)$ by

$$\mathbf{R}(s,t) := \sum_{n=1}^{\infty} \mathbf{K}^{*n}(s,t)$$

(4.16)

where $\mathbf{K}$ is the kernel defined in (4.7).
Then the solution $x(.)$ of (4.6, 4.7, 4.8) is given by

$$x(s) = \mathbf{h}(s) + \int_0^s \mathbf{R}(s,t)\mathbf{h}(t)dt$$

(4.17)

Furthermore, the resolvent kernel $\mathbf{R}$ satisfies

$$\mathbf{R} * \mathbf{K} = \mathbf{K} * \mathbf{R} = \mathbf{R} - \mathbf{K}$$

(4.18)

The <u>proof</u> is analogous to he proof for the case of continuous kernels. ///

We turn to the issue of duality. For ordinary (non-impulsive) linear Volterra equations, duality says that if the equation

$$x(s) = h(s) + \int_0^s K(s,t)x(t)dt$$

has resolvent kernel $R(s,t)$, so that, for every function $h(.)$, the solution is

$$x(s) = h(s) + \int_0^s R(s,t)h(t)dt$$

then the solution of the dual equation

$$y(t) = h(t) + \int_t^T y(s)K(s,t)ds$$

is given, for every $h(.)$, by

$$y(t) = h(t) + \int_t^T h(s)R(s,t)ds$$

For applications to optimal control and variational problems, we need to establish the corresponding results for impulsive Volterra equations.

We have

<u>Theorem 4.3.</u> If the resolvent kernel of equation (4.6) is **R**(s,t), so that, for every h(.) in C([0, T] $\mapsto$ IR), the solution of

$$x(s) = h(s) + \int_0^s \mathbf{K}(s,t) x(t) dt \tag{4.19}$$

is given by

$$x(s) = h(s) + \int_0^s \mathbf{R}(s,t) h(t) dt \tag{4.20}$$

then, for every $\eta$ in C(0,T;$\boldsymbol{\tau}$), the function

$$y(t) := \eta(t) + \int_t^T \eta(s) \mathbf{R}(s,t) ds \tag{4.21}$$

satisfies the integral equation

$$y(t) = \eta(t) + \int_t^T y(s) K(s,t) ds + \\ + \sum_{j: t < \tau_j < T} \sum_{i=j}^N \sum_{\sigma \in \mathbf{P}(j,i)} K(\tau_j, t) \Lambda(\sigma) \int_t^T y(s) L(s, \tau_i) ds \tag{4.22}$$

<u>Proof:</u> The function y given by (4.20) satisfies

$$y(t) = \eta(t) + \int_t^T y(s)\mathbf{K}(s,t)ds$$

(4.23)

This follows from the properties of **R** stated in theorem 4.2, by a calculation analogous to the case of continuous kernels. In view of the representation (4.7) of the kernel **K**, and by interchanging the summation over i and j, i.e., for s>t,

$$\sum_{i,j,\sigma} L(s,\tau_i)\Lambda(\sigma)K(\tau_j,t) = \sum_{j:t<\tau_j<T} \sum_{i=j}^{N} \sum_{\sigma\in\mathbf{P}(j,i)} K(\tau_j,t)\Lambda(\sigma)L(s,\tau_i)$$

(4.24)

it follows that (4.23) is tantamount to (4.22). ////

## 5. Variational equations.

We denote by $a_i$ the admissible variations in the controls $a_i$. Under a change of the control $u(t)$ from the values $a_i$ to $a_i+\varepsilon\alpha_i$, the state $x(t)$ changes to $x_\varepsilon(t)$, and the cost functional $J$ changes to $J_\varepsilon$. We are interested in calculating the variations

$$\delta x(t) := \lim_{\varepsilon \to 0^+} \frac{1}{\varepsilon}(x_\varepsilon(t) - x(t)), \quad \delta J := \lim_{\varepsilon \to 0^+} \frac{1}{\varepsilon}(J_\varepsilon - J)$$

(5.1)

The equations for the variations $\delta x$ and $\delta J$ are found by standard methods of the calculus of variations:

$$\delta x(s) = \int_0^s \{f_x(s,t,x(t),u(t))\delta x(t) + f_u(s,t,x(t),u(t))\delta u(t)\}dt +$$

$$+ \sum_{i:0<\tau_i<s} \{G_x(s,\tau_i,x(\tau_i^-),u(\tau_i^-),u(\tau_i^+))\delta x(\tau_i^-) +$$

$$+ G_a(s,\tau_i,x(\tau_i^-),u(\tau_i^-),u(\tau_i^+))\alpha_i + G_b(s,\tau_i,x(\tau_i^-),u(\tau_i^-),u(\tau_i^+))\alpha_{i+1}\}$$

(5.2)

where

$$\delta u(t) = \alpha_i \text{ for } \tau_{i-1} \leq t < \tau_i;$$

$$G_a(s,\tau,x,a,b) = \frac{\partial G(s,\tau,x,a,b)}{\partial a}, \quad G_b(s,\tau,x,a,b) = \frac{\partial G(s,\tau,x,a,b)}{\partial b}$$

(5.3)

For fixed $x(.)$ and $u(.)$, Eq. (5.2) is a linear impulsive Volterra equation, and it is convenient to write it in the form

$$\delta x(s) = \eta(s) + \int_0^s f_x(s,t,x(t),u(t))\delta x(t)dt + \sum_{i:0<\tau_i<s} G_x(s,\tau_i,\ldots)\delta x(\tau_i^-)$$

(5.4)

where

$$\eta(s) = \int_0^s f_u(s,t,x(t),u(t))\delta u(t)dt +$$

$$+ \sum_{i:0<\tau_i<s} \{G_a(s,\tau_i,x(\tau_i^-),u(\tau_i^-),u(\tau_i^+))\alpha_i + G_b(s,\tau_i,x(\tau_i^-),u(\tau_i^-),u(\tau_i^+))\alpha_{i+1}\}$$

(5.5)

We note that

$$\int_0^s f_u(s,t,x(t),u(t))\delta u(t)dt = \sum_{i=1}^{N_s} \int_{\tau_{i-1}}^{\tau_i} f_u(s,t,x(t),a_i)\alpha_i dt +$$

$$+ \int_{\tau_{N_s}}^s f_u(s,t,x(t),a_{N_s+1})\alpha_{N_s+1}dt$$

(5.6)

The variation of J is

$$\delta J = \int_0^T \{g_x(t,x(t),u(t))\delta x(t) + g_u(t,x(t),u(t))\delta u(t)\}dt +$$

$$+ \sum_{i=1}^{N+1} \{\frac{\partial \varphi(x(\tau_1^-),x(\tau_2^-),...,x(\tau_{N+1}^-),a_1,a_2,...,a_{N+1})}{\partial x_i}\delta x(\tau_i^-) +$$

$$+ \frac{\partial \varphi(x(\tau_1^-),x(\tau_2^-),...,x(\tau_{N+1}^-),a_1,a_2,...,a_{N+1})}{\partial a_i}\alpha_i\}$$

(5.7)

and we have

$$\int_0^T g_u(t,x(t),u(t))\delta u(t)dt = \sum_{i=1}^{N+1} \int_{\tau_{i-1}}^{\tau_i} g_u(t,x(t),a_i)\alpha_i dt$$

(5.8)

The variational equation (5.4) is of the type (4.1) with

$$K(s,t) = f_x(s,t,x(t),u(t)), \quad L(s,\tau_i) = G_x(s,\tau_i,x(\tau_i^-),a_i,a_{i+1})$$

(5.9)

and the solution dx(t) can be represented in terms of the resolvent R defined in section 4 of this paper,

$$\delta x(s) = \eta(s) + \int_0^s \mathbf{R}(s,t)\eta(t)dt$$

(5.10)

Acording to the results of section 4, (5.10) implies

$$\delta x(\tau_i^-) = \sum_{j=1}^{i} \sum_{\sigma \in \mathbf{P}(j,i)} \Gamma(\sigma)[\eta(\tau_j^-) + \int_0^{\tau_j} K(\tau_j,t)\delta x(t)dt]$$

(5.11)

or equivalently, in view of (5.5) and (5.9),

$$\delta x(\tau_i^-) = \sum_{j=1}^{i} \sum_{\sigma \in P(j,i)} \Gamma(\sigma)[\int_0^{\tau_j} \{f_x(\tau_j,t,x(t),u(t))\delta x(t) + f_u(\tau_j,t,x(t),u(t))\delta u(t)\}dt +$$
$$+ \sum_{k=1}^{j-1} \{G_a(s,\tau_k,x(\tau_k^-),a_k,a_{k+1})\alpha_k + G_b(s,\tau_k,x(\tau_k^-),a_k,a_{k+1})\alpha_{k+1}\}]$$

(5.12)

The resolvent kernel **R** is obtained from the continuous and impulsive kernels given by (5.9). The discrete resolvent kernel $\Gamma(\sigma)$ corresponding to any increasing path $\sigma := (\tau_j, \tau_{k_1}, \tau_{k_2}, ..., \tau_{k_\alpha}, \tau_i)$ is given by

$$\Gamma(\sigma) = G_x(\tau_i, \tau_{k_\alpha}, x(\tau_{k_\alpha}^-), a_{k_\alpha}, a_{k_\alpha+1})G_x(\tau_{k_\alpha}, \tau_{k_{\alpha-1}}, x(\tau_{k_{\alpha-1}}^-), a_{k_{\alpha-1}}, a_{k_{\alpha-1}+1}) \bullet$$
$$\bullet \cdots \bullet G_x(\tau_{k_1}, \tau_j, x(\tau_j^-), a_j, a_{j+1})$$

(5.13)

By substituting (5.12) into (5.7), we obtain

$$\delta J = \int_0^T [g_x(t,x(t),u(t)) +$$

$$+ \sum_{i=1}^{N+1} \sum_{j=1}^{i} \sum_{\sigma \in \mathbf{P}(j,i)} \varphi_{x_i}(x(\tau_1^-),\ldots,x(\tau_{N+1}^-),a_1,\ldots,a_{N+1})\Gamma(\sigma)f_x(\tau_j,t,x(t),u(t))]\delta x(t)dt +$$

$$+ \int_0^T g_u(t,x(t),u(t))\delta u(t)dt + \sum_{i=1}^{N+1} \varphi_{a_i}(x(\tau_1^-),\ldots,x(\tau_{N+1}^-),a_1,\ldots,a_{N+1})\alpha_i +$$

$$+ \sum_{i=1}^{N+1} \sum_{j=1}^{i} \sum_{\sigma \in \mathbf{P}(j,i)} \varphi_{x_i}(x(\tau_1^-),\ldots,x(\tau_{N+1}^-),a_1,\ldots,a_{N+1})\Gamma(\sigma)[\int_0^{\tau_j} f_u(\tau_j,t,x(t),u(t))\delta u(t)dt +$$

$$+ \sum_{k=1}^{j-1} \{G_a(s,\tau_k,x(\tau_k^-),a_k,a_{k+1})\alpha_k + G_b(s,\tau_k,x(\tau_k^-),a_k,a_{k+1})\alpha_{k+1}\}]$$

(5.14)

For later reference, we also write down the explicit form of (5.10), namely

$$\delta x(s) = \int_0^s f_u(s,t,x(t),u(t))\delta u(t)dt +$$

$$+ \sum_{i:0<\tau_i<s} \{G_a(s,\tau_i,x(\tau_i^-),a_i,a_{i+1})\alpha_i + G_b(s,\tau_i,x(\tau_i^-),a_i,a_{i+1})\alpha_{i+1}\} +$$

$$+ \int_{r=0}^s \int_{t=0}^r \mathbf{R}(s,r)f_u(r,t,x(t),u(t))\delta u(t)dtdr +$$

$$+ \int_{r=0}^s \sum_{i:0<\tau_i<r} \mathbf{R}(s,r)\{G_a(r,\tau_i,x(\tau_i^-),a_i,a_{i+1})\alpha_i + G_b(r,\tau_i,x(\tau_i^-),a_i,a_{i+1})\alpha_{i+1}\}dr$$

(5.15)

where the double integral in (5.15) can also be written as

$$\int_{t=0}^s \int_{r=t}^s \mathbf{R}(s,r)f_u(r,t,x(t),u(t))\delta u(t)drdt \ .$$

## 6. Maximum principle.

We define

$$\xi(t) := g_x(t, x(t), u(t)) +$$
$$+ \sum_{i=1}^{N+1} \sum_{j=1}^{i} \sum_{\sigma \in \mathbf{P}(j,i)} f_x(\tau_j, t, x(t), u(t)) \Gamma(\sigma) \varphi_{x_i}(x(\tau_1^-), ..., x(\tau_{N+1}^-), a_1, ..., a_{N+1})$$

(6.1)

$$\psi(t) = -\xi(t) - \int_t^T \xi(s) \mathbf{R}(s,t) ds$$

(6.2)

Then we have, by the duality results of section 4,

$$\psi(t) = -g_x(t, x(t), u(t)) -$$
$$- \sum_{i=1}^{N+1} \sum_{j=1}^{i} \sum_{\sigma \in \mathbf{P}(j,i)} f_x(\tau_j, t, x(t), u(t)) \Gamma(\sigma) \varphi_{x_i}(x(\tau_1^-), ..., x(\tau_N^-), a_1, ..., a_{N+1}) +$$
$$+ \int_t^T \psi(s) f_x(s, t, x(t), u(t)) dt +$$
$$+ \sum_{j: t < \tau_j < T} \sum_{i=j}^{N+1} \sum_{\sigma \in \mathbf{P}(j,i)} f_x(\tau_j, t, x(t), u(t)) \Gamma(\sigma) \int_t^T \psi(s) G_x(s, \tau_i, x(\tau_i^-), u(\tau_i^-), u(\tau_i^+)) ds$$

(6.3)

The expression for $\delta J$ can be written as

$$\delta J = \int_0^T \xi(t)\delta x(t)dt + \int_0^T g_u(t,x(t),u(t))\delta u(t)dt +$$

$$+ \sum_{i=1}^{N+1} \varphi_{a_i}(x(\tau_1^-),...,x(\tau_{N+1}^-),a_1,...,a_{N+1})\alpha_i +$$

$$+ \sum_{i=1}^{N+1} \sum_{j=1}^{i} \sum_{\sigma \in \mathbf{P}(j,i)} \varphi_{x_i}(x(\tau_1^-),...,x(\tau_{N+1}^-),a_1,...,a_{N+1})\Gamma(\sigma)\eta(\tau_j^-)$$

(6.4)

We use $\delta x(t) = \eta(t) + \int_0^t \mathbf{R}(t,t_1)\eta(t_1)dt_1$ in (6.4) and we find

$$\delta J = \int_0^T \xi(t)\eta(t)dt + \int_0^T \int_t^T \xi(t_1)\mathbf{R}(t_1,t)dt_1 \eta(t)dt +$$

$$+ \int_0^T g_u(t,x(t),u(t))\delta u(t)dt +$$

$$+ \sum_{i=1}^{N+1} \varphi_{a_i}(x(\tau_1^-),...,x(\tau_N^-),a_1,...,a_{N+1})\alpha_i +$$

$$+ \sum_{i=1}^{N+1} \sum_{j=1}^{i} \sum_{\sigma \in \mathbf{P}(j,i)} \varphi_{x_i}(x(\tau_1^-),...,x(\tau_N^-),a_1,...,a_{N+1})\Gamma(\sigma)\eta(\tau_j^-) =$$

$$= -\int_0^T \psi(t)\eta(t)dt + \int_0^T g_u(t,x(t),u(t))\delta u(t)dt +$$

$$+ \sum_{i=1}^{N+1} \varphi_{a_i}(x(\tau_1^-),...,x(\tau_N^-),a_1,...,a_{N+1})\alpha_i +$$

$$+ \sum_{i=1}^{N+1} \sum_{j=1}^{i} \sum_{\sigma \in \mathbf{P}(j,i)} \varphi_{x_i}(x(\tau_1^-),...,x(\tau_N^-),a_1,...,a_{N+1})\Gamma(\sigma)\eta(\tau_j^-)$$

(6.5)

In view of the results of section 4, we have

$$\int_0^T \psi(t)\eta(t)dt = \int_0^T \int_t^T \psi(s)f_u(s,t,x(t),u(t))ds\delta u(t)dt +$$

$$+ \int_0^T \psi(t) \sum_{i=1}^N \{G_a(t,\tau_i,x(\tau_i^-),u(\tau_i^-),u(\tau_i^+))\alpha_i +$$

$$+ G_b(t,\tau_i,x(\tau_i^-),u(\tau_i^-),u(\tau_i^+))\alpha_{i+1}\}dt =$$

$$= \sum_{i=1}^{N+1} [\int_{\tau_{i-1}}^{\tau_i} \alpha_i \int_t^T \psi(s)f_u(s,t,x(t),u(t))dsdt +$$

$$+ \int_{\tau_i}^T \psi(t)\{\alpha_i G_a(t,\tau_i,x(\tau_i^-),u(\tau_i^-),u(\tau_i^+)) +$$

$$+ \alpha_{i+1} G_b(t,\tau_i,x(\tau_i^-),u(\tau_i^-),u(\tau_i^+))\}dt]$$

(6.6)

Further, we have (cf. (5.8))

$$\int_0^T g_u(t,x(t),u(t))\delta u(t)dt = \sum_{i=1}^{N+1} \alpha_i \int_{\tau_{i-1}}^{\tau_i} g_u(t,x(t),a_i)dt;$$

$$\sum_{i=1}^{N+1} \sum_{j=1}^{i} \sum_{\sigma \in \mathbf{P}(j,i)} \varphi_{x_i}(x(\tau_1^-),...,x(\tau_{N+1}^-),a_1,...,a_{N+1})\Gamma(\sigma)\eta(\tau_j^-) =$$

$$= \sum_{i=1}^{N+1} \sum_{j=1}^{i} \sum_{\sigma \in \mathbf{P}(j,i)} \varphi_{x_k}(x(\tau_1^-),...,x(\tau_{N+1}^-),a_1,...,a_{N+1})\Gamma(\sigma) \cdot$$

$$\cdot [\sum_{k=1}^{j} \alpha_k \int_{\tau_{k-1}}^{\tau_k} f_u(\tau_j,t,x(t),u(t))dt +$$

$$+ \sum_{k=1}^{j-1} \{\alpha_k G_a(\tau_j,\tau_k,x(\tau_k^-),u(\tau_k^-),u(\tau_k^+)) +$$

$$+ \alpha_{k+1} G_b(\tau_j,\tau_k,x(\tau_k^-),u(\tau_k^-),u(\tau_k^+))\}]$$

(6.7)

Thus, we have:

$$\delta J = -\sum_{i=1}^{N+1} \Delta_i \alpha_i \quad (6.8)$$

where the terms $\Delta_i$ are expressed in terms of x, u, $\psi$, G, $\varphi$, as follows:

$$\Delta_i = \int_{\tau_{i-1}}^{\tau_i} \int_t^T \psi(s) f_u(s, t, x(t), u(t)) ds\, dt +$$

$$+ \int_{\tau_i}^T \psi(t) G_a(t, \tau_i, x(\tau_i^-), u(\tau_i^-), u(\tau_i^+)) dt +$$

$$+ \int_{\tau_{i-1}}^T \psi(t) G_b(t, \tau_{i-1}, x(\tau_{i-1}^-), u(\tau_{i-1}^-), u(\tau_{i-1}^+)) dt -$$

$$- \varphi_{a_i}(x(\tau_1^-), ..., x(\tau_{N+1}^-), a_1, ..., a_{N+1}) - \int_{\tau_{i-1}}^{\tau_i} g_u(t, x(t), a_i) dt -$$

$$- \sum_{k=i}^{N+1} \sum_{j=i}^{k} \sum_{\sigma \in \mathbf{P}(j,k)} \varphi_{x_i}(x(\tau_1^-), ..., x(\tau_{N+1}^-), a_1, ..., a_{N+1}) \Gamma(\sigma) \cdot$$

$$\cdot [\int_{\tau_{i-1}}^{\tau_i} f_u(\tau_j, t, x(t), u(t)) dt + G_b(\tau_j, \tau_{i-1}, x(\tau_{i-1}^-), u(\tau_{i-1}^-), u(\tau_{i-1}^+))] -$$

$$- \sum_{k=i-1}^{N+1} \sum_{j=i-1}^{k} \sum_{\sigma \in \mathbf{P}(j,k)} \varphi_{x_i}(x(\tau_1^-), ..., x(\tau_{N+1}^-), a_1, ..., a_{N+1}) \Gamma(\sigma) \cdot$$

$$\cdot G_a(\tau_j, x(\tau_i^-), u(\tau_i^-), u(\tau_i^+)) \quad (6.9)$$

The necessary condition for a minimum of J, namely

$\delta J \geq 0$ for every admissible $\delta u$

is tantamount to

$$\sum_{i=1}^{N+1} \Delta_i \alpha_i \leq 0 \text{ for every admissible } \boldsymbol{\alpha} = (\alpha_1, \ldots, \alpha_{N+1})$$

(6.10)

We note that the expression

$$H(t, x, u, \psi(\cdot)) := -g(t, x(t), u(t)) -$$
$$- \sum_{i=1}^{N+1} \sum_{j=1}^{i} \sum_{\sigma \in \mathbf{P}(j,i)} f(\tau_j, t, x(t), u(t)) \Gamma(\sigma) \varphi_{x_i}(x(\tau_1^-), \ldots, x(\tau_{N+1}^-), a_1, \ldots, a_{N+1}) +$$
$$+ \int_t^T \psi(s) f(s, t, x(t), u(t)) ds +$$
$$+ \sum_{j: t < \tau_j < T} \sum_{i=j}^{N+1} \sum_{\sigma \in \mathbf{P}(j,i)} f(\tau_j, t, x(t), u(t)) \Gamma(\sigma) \int_t^T \psi(s) G_x(s, \tau_i, x(\tau_i^-), u(\tau_i^-), u(\tau_i^+)) ds$$

(6.11)

is the impulsive Hamiltonian, in the sense that the co-state $\psi$ satisfies

$$\psi(t) = \frac{\partial}{\partial x} H(t, x, u, \psi(\cdot))$$

(6.12)

where, in the partial differentiation operator $\frac{\partial}{\partial x}$, the symbol "x" refers to the slots, in the expression of H, that are occupied by x(t), but not to the terms $x(\tau_i^-)$.

On the basis of the results above, we have

<u>Theorem 6.1.</u> (The maximum principle for impulsive Volterra equations with discrete controls.) Under the conditions of continuous differentiability of the functions f, G, g, φ with respect to x, u, **a**, assuming the existence of an optimal control policy, we denote by **a**$^*$ an optimal choice for **a**, by $u^*(\cdot)$ the associated control law as a function of t, and by $x^*(\cdot)$ the associated solution of the Volterra integral equation (2.). We denote by $\psi^*(\cdot)$ the solution of

$$\psi^*(t) = \frac{\partial}{\partial x} H(t, x^*, u^*, \psi^*(\cdot))$$

where H is defined by (6.11). Then the first-order necessary condition for optimality is

$$\sum_{i=1}^{N+1} \Delta_i^* \alpha_i \geq 0 \text{ for all admissible } \alpha_i,$$

where $\Delta_i^*$ is calculated from (6.) with $\mathbf{a}^*$, $u^*(\cdot)$, $x^*(\cdot)$, and $\psi^*(\cdot)$ in lieu of $\mathbf{a}$, u, x, and $\psi$, respectively. ///

## 7. Application to impulsive differential equations.

In this section, we examine how the results of section 6 can be specialized, from the general case of impulsive Volterra equations to the particular case of impulsive ordinary differential equations. Impulsive control problems for systems governed by ordinary differential equations have been studied in [M] and the further references therein. The maximum principle for impulsive ODE control problems is contained in [M], but in a form different from the results we state in this section. The problem we consider here has the impulses, in both the state dynamics and the cost functional, dependent on the current and the switched value of the control (i.e. the functions G and $\Phi$ depend on both $a_i$ and $a_{i+1}$); this type of controls amounts to having the admissible values of the control $w_i \equiv (a_i, a_{i+1})$ dependent on the previous value $w_{i-1}$, a type of constraint that does not fall within the scope of [M] and other references on impulsive control of ODE. Our contribution in this section is to demonstrate that the maximum principle (conditions (7.13) below) for impulsive ODE with impulsive controls can be obtained as a particular case of the maximum principle for controlled impulsive Volterra equations proved in the previous section.

Now, we consider the controlled impulsive ODE

$$\frac{dx(t)}{dt} = f(t, x(t), u(t)), \text{ for } t \neq \tau_i; \ x(0) = x_0;$$

$$x(\tau_i^+) = x(\tau_i^-) + G(\tau_i, x(\tau_i^-), a_i, a_{i+1}), \ i = 1, 2, ..., N$$

(7.1)

which can also be written in integral form

$$x(t) = x_0 + \int_0^t f(s, x(s), u(s))ds + \sum_{i: 0 < \tau_i < t} G(\tau_i, x(\tau_i^-), a_i, a_{i+1})$$

(7.2)

We shall use the notation and terminology of sections 2-6 above, with the indicated specializations. The cost functional J is taken in the form

$$J := \int_0^T g(t, x(t), u(t))dt + \sum_{i=1}^N \Phi(x(\tau_i^-), a_i, a_{i+1}) + g_0(x(T))$$

(7.3)

The functional J is plainly a particular case of the functional (2.4), with

$$\varphi(x(\tau_1^-),\ldots,x(\tau_{N+1}^-),a_1,\ldots,a_{N+1}) = \sum_{i=1}^{N} \Phi(x(\tau_i^-),a_i,a_{i+1}) + g_0(x(\tau_{N+1}^-))$$

(7.4)

The co-state $p(t)$ for the problem (7.2-7.3) is defined in terms of the co-state $\psi(t)$ introduced in section 7:

$$p(t) = \int_t^T \psi(s)\,ds + \sum_{j:t<\tau_j<T}\sum_{i=j}^{N}\sum_{\sigma\in \mathbf{P}(j,i)} \Gamma(\sigma)\cdot$$
$$\cdot[-\Phi_{x_i}(x(\tau_i^-),a_i,a_{i+1}) + G_x(\tau_i,x(\tau_i^-),a_i)\int_{\tau_i}^T \psi(s)\,ds] -$$
$$-\sum_{j:t<\tau_j\leq T}\sum_{\sigma\in \mathbf{P}(j,N+1)} \Gamma(\sigma) g_{0,x}(x(\tau_{N+1}^-))$$

(7.5)

We define the continuous and the discrete Hamiltonians, $H_c$ and $H_d$, respectively, for this problem, by

$$H_c(t,x,p,a) := g(t,x,a) - p\,f(t,x,a);$$
$$H_d(t,x,p,a,b) := \Phi(t,x,a,b) - p\,G(t,x,a,b)$$

(7.6)

We shall show that the co-state $p(t)$, defined by (7.5), satisfies

$$\frac{dp(t)}{dt} = \frac{\partial H_c(t, x(t), p(t), u(t))}{\partial x}, \text{ for } t \neq \tau_i, \ i = 0,1,2,..., N+1;$$

$$p(\tau_k^+) - p(\tau_k^-) = \frac{\partial H_d(\tau_k, x(\tau_k^-), p(\tau_k^+), a_k, a_{k+1})}{\partial x} \text{ for } k = 1,2,..., N;$$

$$p(\tau_{N+1}^-) = g_{0,x}(x(\tau_{N+1}^-))$$

(7.7)

First, we establish the jump condition on p, i.e. the second equation in (7.7). It follows from the definition of p that

$$p(\tau_k^+) - p(\tau_k^-) = -\sum_{i=k}^{N} \sum_{\sigma \in \mathbf{P}(k,i)} \Gamma(\sigma) \cdot$$

$$\cdot [-\Phi_{x_i}(\tau_i, x(\tau_i^-), a_i, a_{i+1}) + G_x(\tau_i, x(\tau_i^-), a_i, a_{i+1}) \int_{\tau_i}^{T} \psi(s)ds] +$$

$$+ \sum_{\sigma \in \mathbf{P}(k, N+1)} \Gamma(\sigma) g_{0,x}(x(\tau_{N+1}^-))$$

(7.8)

and

$$p(\tau_k^+) = \int_{\tau_k}^{T} \psi(s)ds +$$

$$+ \sum_{j=k+1}^{N} \sum_{i=j}^{N} \sum_{\sigma \in \mathbf{P}(j,i)} \Gamma(\sigma)[-\Phi_{x_i}(\tau_i, x(\tau_i^-), a_i, a_{i+1}) + G_x(\tau_i, x(\tau_i^-), a_i, a_{i+1}) \int_{\tau_i}^{T} \psi(s)ds] -$$

$$- \sum_{j=k+1}^{N+1} \sum_{\sigma \in \mathbf{P}(j, N+1)} \Gamma(\sigma) g_{0,x}(x(\tau_{N+1}^-))$$

(7.9)

Thus, in order to verify that

$$p(\tau_k^+) - p(\tau_k^-) = \Phi_{x_k}(\tau_k, x(\tau_k^-), a_k, a_{k+1}) - p(\tau_k^+) G_x(\tau_k, x(\tau_k^-), a_k, a_{k+1})$$

(7.10)

it suffices to prove that

$$\sum_{i=k}^{N} \sum_{\sigma \in \mathbf{P}(k,i)} \Gamma(\sigma) \cdot$$

$$\cdot [\Phi_{x_i}(\tau_i, x(\tau_i^-), a_i, a_{i+1}) - G_x(\tau_i, x(\tau_i^-), a_i, a_{i+1}) \int_{\tau_i}^{T} \psi(s)ds] +$$

$$+ \sum_{\sigma \in \mathbf{P}(k,N+1)} \Gamma(\sigma) g_{0,x}(x(\tau_{N+1}^-))\} =$$

$$= \Phi_{x_k}(x(\tau_k^-), a_k, a_{k+1}) - G_x(\tau_k, x(\tau_k^-), a_k, a_{k+1}) \cdot$$

$$\cdot \{\int_{\tau_k}^{T} \psi(s)ds +$$

$$+ \sum_{j=k+1}^{N} \sum_{i=j}^{N} \sum_{\sigma \in \mathbf{P}(j,i)} \Gamma(\sigma)[-\Phi_{x_i}(\tau_i, x(\tau_i^-), a_i, a_{i+1}) + G_x(\tau_i, x(\tau_i^-), a_i, a_{i+1}) \int_{\tau_i}^{T} \psi(s)ds] -$$

$$- \sum_{j=k+1}^{N+1} \sum_{\sigma \in \mathbf{P}(j,N+1)} \Gamma(\sigma) g_{0,x}(x(\tau_{N+1}^-))\}$$

(7.11)

Now, every $\sigma \in \mathbf{P}(k,i)$ is one of the following: if i=k, then $\sigma = \{\tau_k\}$; if i>k, then $\sigma$ is either $\{\tau_k, \tau_i\}$ or the concatenation of $\{\tau_k, \tau_j\}$, for some j, $k < j \leq i$, with some $\sigma' \in \mathbf{P}(j,i)$. Correspondingly, we have:

if $\sigma = \{\tau_k\}$, then $\Gamma(\sigma) = 1$, by definition;

if $\sigma = \{\tau_k, \tau_i\}$, then $\Gamma(\sigma) = G_x(\tau_k, x(\tau_k^-), a_k, a_{k+1})$;

if $\sigma$ is the concatenation of $\{\tau_k, \tau_j\}$, for some j, $k < j \leq i$, with some $\sigma' \in \mathbf{P}(j,i)$, then

$\Gamma(\sigma) = G_x(\tau_k, x(\tau_k^-), a_k, a_{k+1}) \Gamma(\sigma')$.

Therefore,

$$\sum_{i=k}^{N} \sum_{\sigma \in \mathbf{P}(k,i)} \Gamma(\sigma) \cdot$$

$$\cdot [\Phi_{x_i}(\tau_i, x(\tau_i^-), a_i, a_{i+1}) - G_x(\tau_i, x(\tau_i^-), a_i, a_{i+1}) \int_{\tau_i}^{T} \psi(s) ds] +$$

$$+ \sum_{j=k}^{N+1} \sum_{\sigma \in \mathbf{P}(j,N+1)} \Gamma(\sigma) g_{0,x}(x(\tau_{N+1}^-)) =$$

$$= \Phi_{x_k}(x(\tau_k^-), a_k, a_{k+1}) - G_x(\tau_k, x(\tau_k^-), a_k, a_{k+1}) \cdot \int_{\tau_k}^{T} \psi(s) ds +$$

$$+ \sum_{i=k+1}^{N} \sum_{j'=k+1}^{i} \sum_{\sigma' \in \mathbf{P}(j',i)} \Gamma(\sigma') G_x(\tau_k, x(\tau_k^-), a_k, a_{k+1}) \cdot$$

$$\cdot [\Phi_{x_i}(\tau_i, x(\tau_i^-), a_i, a_{i+1}) - G_x(\tau_i, x(\tau_i^-), a_i, a_{i+1}) \int_{\tau_i}^{T} \psi(s) ds] -$$

$$- \sum_{j'=k+1}^{N+1} \sum_{\sigma' \in \mathbf{P}(j',N+1)} G_x(\tau_k, x(\tau_k^-), a_k, a_{k+1}) \cdot$$

$$\cdot \Gamma(\sigma') g_{0,x}(x(\tau_{N+1}^-)) \}$$

(7.12)

and then (7.11) follows by a change in the order of summation, namely

$$\sum_{i=k+1}^{N} \sum_{j'=k+1}^{i} = \sum_{j'=k+1}^{N} \sum_{i=j'}^{N} .$$

The first equation in (7.7) follows from (6.12) and the observation that, under the stated definitions, for $t \notin \boldsymbol{\tau}$, we have

$$\frac{dp(t)}{dt} = -\psi(t) \quad \text{and} \quad H_c(t, x(t), p(t), u(t)) = -H(t, x, u, \psi(\cdot))$$

where $H(t, x, u, \psi(\cdot))$ is defined in (6.11).
Finally, the condition

$$p(\tau_{N+1}^-) = g_{0,x}(x(\tau_{N+1}^-))$$

follows from the definition of p in (7.5) and the observation that, for k=N+1, the only element of $\mathbf{P}(k, N+1)$ is $\sigma = \{ \tau_{N+1} \}$ for which $\Gamma(\sigma) = 1$ by definition.

Consequently, we have shown:

Theorem 7.1. (Maximum principle for impulsive ordinary differential equations.) Under conditions of continuous differentiability of f, G, g, $\Phi$, G with respect to x, u, **a**, and with the superscript "*" denoting optimality, we have the first-order necessary conditions for the problem (7.1, 7.2) as follows:

$$\frac{dp^*(t)}{dt} = \frac{\partial H_c(t, x^*(t), p^*(t), u^*(t))}{\partial x}, \text{ for } t \neq \tau_i,\ i = 0, 1, 2, \dots, N+1;$$

$$p^*(\tau_k^+) - p^*(\tau_k^-) = \frac{\partial H_d(\tau_k, x^*(\tau_k^-), p^*(\tau_k^+), a_k^*, a_{k+1}^*)}{\partial x} \text{ for } k = 1, 2, \dots, N;$$

$$p^*(\tau_{N+1}^-) = g_{0,x}(x^*(\tau_{N+1}^-));$$

$$\sum_{i=1}^{N+1} \Delta_i^* \alpha_i^* \geq 0 \text{ for all admissible } \alpha_i$$

(7.13)

where each $\Delta_i$ is evaluated from (6.9) with $\psi(t) := -\frac{dp(t)}{dt}$. ///

The impulsive differential equation (7.1) can also be written in the Hamiltonian form

$$\frac{dx^*(t)}{dt} = -\frac{\partial H_c(t, x^*(t), p^*(t), u^*(t))}{\partial p} \text{ for } t \neq \tau_i,\ i = 0, 1, 2, \dots, N+1;$$

$$x^*(\tau_k^+) - x^*(\tau_k^-) = -\frac{\partial H_d(\tau_k, x^*(\tau_k^-), p^*(\tau_k^+), a_k^*, a_{k+1}^*)}{\partial p} \text{ for } k = 1, 2, \dots, N;$$

$$x^*(\tau_0^+) = x_0$$

(7.14)

References

[B].  S. A. BELBAS, Iterative schemes for optimal control of Volterra integral equations, Nonlinear Analysis, Vol. 37, 1999, pp. 57-79.

[BA].  V. L. BAKKE, A maximum principle for an optimal control problem with integral constraints, Journal of Optimization Theory and Applications, Vol. 13, 1974, pp. 32-55.

[BI].  L. BITTNER, Begründung des sogenanten diskreten Maximumprinzips, Z. Wahrscheinlichkeitstheorie und verw. Gebiete, Vol. 10, 1968, pp. 289-304.

[BO].  V. G. BOLTYANSKII, Optimal control of discrete systems, J. Wiley, New York, 1978.

[CA].  D. A. CARLSON, An elementary proof of the maximum principle for optimal control problems governed by a Volterra integral equation, Journal of Optimization Theory and Applications, Vol. 54, 1987, pp. 43-61.

[C].  C. CORDUNEANU, Integral equations and applications, Cambridge U. Press, Cambridge, 1991.

[DM].  B. P. DEMIDOVICH, I. A. MARON, Computational mathematics, Mir Publishers, Moscow, 1987.

[KM].  M. I. KAMIEN, E. MULLER, Optimal control with integral state equations, The Review of Economic Studies, Vol. 43, 1976, pp. 469-473.

[K].  J. KONDO, Integral equations, Kodansa, Tokyo, and Clarendon Press, Oxford, 1991.

[M].  B. M. MILLER, E. YA. RUBINOVICH, Impulsive control in continuous and discrete-continuous systems, Kluwer, New York, 2003.

[RS].  P. RUTKAUSKAS, W. H. SCHMIDT, The optimality conditions for integral processes with delays, Lietuvos Matematikos Rinkinys, Vol. 23, No. 3, 1983, pp. 117-126.

[S1].  W. H. SCHMIDT, Durch Integralgleichungen beschrienbene optimale Prozesse mit Nebenbedingungen in Banachräumen - notwendige Optimalitätsbedingungen, ZAMM, Vol. 62, 1982, pp. 65-75.

[S2].  W. H. SCHMIDT, Necessary optimality conditions for discrete integral control processes with phase constraints and optimal choice of switchings (in Russian), Lietuvos Matematikos Rinkinys, Vol. 23, No. 3, 1983, pp. 190-195.

[S3]. W. H. SCHMIDT, Necessary optimality conditions for discrete integral processes in Banach spaces, <u>Beiträge zur Analysis</u>, Vol. 16, 1981, pp. 137-145.

[S4]. W. H. SCHMIDT, Maximum principles for processes governed by integral equations in Banach spaces as sufficient optimality conditions, <u>Beiträge zur Analysis</u>, Vol. 17, 1981, pp. 85-93.

[VS]. O. VON STRYK, <u>Numerische Lösung optimaler Steuerungsprobleme: Diskretisierung, Parameteroptimierung und Berechnung der adjugierten Variablen</u>, Fortschritt-Berichte VDI, Reihe 8, Nr. 441, VDI-Verlag, Düsseldorf, 1995.